\input amstex
 \documentstyle{amsppt}
 \magnification=1200
 \vsize=8.5truein
 \hsize= 6.0truein
 \hoffset=0.3truein
 \voffset=0.2truein
  \TagsOnRight

\documentstyle{amsppt}\magnification=\magstep1

\TagsOnRight
\NoBlackBoxes

\baselineskip=12pt
\parindent=10pt

\vsize=8.5truein
\hsize=6.0truein
\hoffset=0.3truein
\voffset=0.2truein

\NoRunningHeads


\nologo
 \NoRunningHeads
 \NoBlackBoxes

 \topmatter
\title Certain 4-manifolds with non-negative sectional curvature
\endtitle
\author Jianguo Cao\footnote{Supported partially by NSF Grant
DMS 0404558 \hfill{$\,$}}
\endauthor
\address Mathematics Department, University of Notre Dame,
Notre Dame, IN 46556
\endaddress
\email jcao\@nd.edu
\endemail

\abstract

     In this paper, we study certain compact 4-manifolds with non-negative sectional curvature $K$. If $s$
     is the scalar curvature and $W_+$ is the self-dual part of Weyl tensor, then it will be shown that there
     is no metric $g$ on $S^2 \times S^2$ with both (i) $K > 0$ and (ii) $ \frac{1}{6} s - W_+ \ge 0$.

        We also investigate other aspects of 4-manifolds with non-negative sectional curvature. One of our results
        implies a theorem of Hamilton: ``If a simply-connected, closed 4-manifold $M^4$ admits a metric $g$ of
        non-negative curvature operator, then $M^4$ is one of $S^4$, $\Bbb CP^2$ and $S^2 \times S^2$".

        Our method is different from Hamilton's and is much simpler. A new version of the second variational formula for
        minimal surfaces in $4$-manifolds is proved.


\endabstract
\endtopmatter

\document
 \baselineskip=.25in

     \vskip2mm

\vskip2mm
\vskip2mm

 \head  Introduction
 \endhead

    \vskip4mm

In this paper, we are interested in 4-dimensional Riemannian manifolds of non-negative sectional curvature. For simplicity,
we denote the sectional curvature by $K$. By $K \ge 0$, we mean that all sectional curvatures are non-negative.

   First, let us say a few words about non-compact $4$-manifolds of non-negative sectional curvature. According to a theorem
   of Cheeger-Gromoll, any complete non-compact $4$-manifold of $K \ge 0$ is diffeomorphic
   to a vector bundle over a compact totally geodesic submanifold
   $N$ of dimension $ \le 3$, (cf. [CG]). R. Hamilton classified
   all $3$-manifolds of non-negative Ricci curvature (cf. [H] Theorem 1.2).  Therefore, {\it  non-compact} complete
   $4$-manifolds of non-negative sectional curvature have been classified and fairly understood.

        In what follows, we will only consider compact simply-connected closed $4$-manifolds of non-negative
        curvature. There is a uniform upper bound for Betti numbers of non-negatively curved manifolds given
by Gromov (cf. [G$r_1$]). This fact and a theorem of Freedman imply that there are only finitely many
non-homeomorphic $4$-manifolds of non-negative sectional curvature (cf. [L] p4-5). It would be interesting to
know which $4$-manifolds can carry a metric of $K \ge 0$.

     We would like to mention some relevant facts which indicates some difficulties to classify all non-negatively
     curved $4$-manifolds. First of all, Sha, Yang and Anderson recently showed that $  \#^k_1 (S^2 \times S^2)$
     and $\#^k_1 ( \Bbb CP^2 \times \Bbb CP^2  )$ admit metrics of strictly positive Ricci curvature for every
     $k \ge 1$, where $M\#N$ means the connected sum of $M$ and $N$ (cf. [SY] Theorem 4). Tian and Yau have
     also proved that $\Bbb CP^2 \# k \overline{\Bbb CP}^2 $ carries K\"ahler-Einstein metrics of
     positive Ricci
     curvature for $ 3 \le k \le 8 $ (cf. [TY]). Furthermore, J. Cheeger in fact construct metrics on
     $\Bbb CP^2 \# \Bbb CP^2$ and $\Bbb CP^2 \# \overline{\Bbb CP}^2$ of
      non-negative sectional curvature (cf. [Ch]).

      Hence, our first attempt is to study non-negatively curved $4$-manifolds with an extra condition.

\proclaim{Theorem A} Let $M$ be a simply-connected closed $4$-manifold whose intersection form is indefinite.
Suppose $M$ carries a smooth metric $g$ of non-negative sectional curvature and
$$
 \quad \frac s6 - W_+ \ge 0, \quad (\text{ or } \quad  \frac s6 - W_- \ge 0). \tag*
$$
Then $M$ is diffeomorphic to $S^2 \times S^2$ and $g$ is a product of metrics of $K \ge 0$.

\endproclaim

   We only consider $4$-manifolds of indefinite type, since any $1$-connected smooth closed $4$-manifold
   $M$ of definite type is homotopy equivalent to $\#^k_1 \Bbb CP^2 $ for some $k \ge 1$ (cf. [L] p3-6).
   Poon showed that $\#^k_1 \Bbb CP^2 $ carries a metric $g$ of $\frac s6 - W_+ \ge 0$ for all $k \le 3$.
   LeBrun recently constructed metrics of $\frac s6 - W_+ \ge 0$ on $\#^k_1 \Bbb CP^2 $ explicitly for
   all $k \ge 1$ (cf. [P], [DF], [L$e_2$]). Hence, any closed $1$-connected smooth $4$-manifold of definite
   type is homotopy equivalent to a Riemannian $4$-manifold $(M^4, g)$ of $\frac s6 - W_+ \ge 0$.

              On $S^2 \times S^2$, one can also construct an $1$-parameter family of twisted metrics which
              satisfy $\frac s6 - W_+ \ge 0$ (cf. section 2).

              Theorem A has an interesting application:

\proclaim{Corollary 0.1} $S^2 \times S^2$ does not carry any metric $g$ of both positive sectional
curvature and $\frac s6 - W_+ \ge 0$.
\endproclaim

     The Hopf problem of whether or not $S^2 \times S^2$ admits a metric $g$ of positive sectional
     curvature is still open. We shall also show

\proclaim{ Theorem B} Let $M$ be a smooth simply-connected closed $4$-manifold which
carries a metric $g$ of non-negative curvature $K \ge 0$ and $\frac s6 - W \ge 0$. Then
one of following holds:

\smallskip
\noindent
(1) $M$ is homeomorphic to $S^4$;

\smallskip
\noindent
(2) $(M, g)$ is complex, K\"ahler and bi-holomorphic to $\Bbb CP^2$,

\smallskip
\noindent
(3) $M$ is diffeomorphic to $S^2 \times S^2$, $g$ is a product of metrics of non-negative
curvature on each factor.
\endproclaim

   Theorem B is useful, since it implies

\proclaim{Corollary 2.2} (Hamilton [H]) Let $M$ be a smooth closed $1$-connected $4$-manifold.
If $M$ carries a metric of non-negative curvature operator. Then one of following holds:

\smallskip
\noindent
(1) M is homeomorphic to $S^4$;

\smallskip
\noindent
(2) $(M, g)$ is complex, K\"ahler and bi-holomorphic to $\Bbb R^n$;

\smallskip
\noindent
(3) $M$ is diffeomorphic to $S^2 \times S^2$, $g$ is a product of metric to $S^2 \times S^2$, $g$ is a product
of metric of non-negative curvature on each factor.

\endproclaim

In case (1) above, Hamilton in fact showed that $M^4$ is diffeomorphic to $S^4$. However, our proof is totally
different from his and is much simpler. Our first observation is following: If $(M^4, g)$ satisfies
$\frac s6 - W_+ \ge 0$ and $M^4$ has the indefinite intersection form on $H_2(M^4, \Bbb R)$, then
$M$ has to be K\"ahler. This can be shown by using the B\"ocher-Weitzenb\"ock formula on $2$-forms. A similar
argument as above was previously used by Bourguignon and LeBrun for a different purpose  (cf. [Bo], [L$e_2$]).
K\"ahler surfaces of non-negative curvature have been classified by Howard-Smyth (cf. [HS]).

     Our second approach is to use the minimal $2$-sphere theory to study the geometry and topology of
     $4$-manifolds. We shall show

\proclaim{Theorem C} Let $S$ be an immersed minimal $2$-sphere in $(S^2 \times S^2, g)$, which is homotopic
to the first factor of $S^2 \times S^2$. Suppose that there is an open set $U$ containing $S$ in which
the metric $g$ satisfies $\frac s6 - W_+ \ge 0$ and $K > 0$, then $S$ is unstable.
\endproclaim

    The existence of some minimal $2$-spheres was given by Sacks and Uhlenbeck. Their result was improved and
    applied to the proof of Frankel's conjecture by Siu-Yau and to the topology of $3$-manifolds by Meeks-Yau
    and Schoen-Yau (cf. [SiY], [ScY]). Some non-existence of minimal $2$-spheres in $4$-manifolds were also
    derived by Futaki (cf. [Fut]). Progress has been made by Micallef, Moore, Wolfson and others
    (cf. [MM], [MW]).

    Although our second approach is less successful than the first one in this paper, the minimal surface
    theory of $4$-manifolds may deserve the further investigation and improvement. Hence, in section 3, we provide
    ourselves with some useful tools. For example, we will derive a new version of the second variational
    formula for minimal surfaces in $4$-manifold. The coefficient in this formula has a term which is
    exactly equal to $\frac s6 - W_+$. Our method is different from [MW]. Some relevant facts may be found in
    [A]  and [Ka].

    \medskip

    {\bf Conventions.} $K$ will always stand for the sectional curvature. We also denote the scalar
    curvature by $s$. $W_+$ (resp. $W_-$) is the self-dual (resp. anti-self-dual) part of the Weyl
    tensor of a given Riemannian $4$-manifold. The curvature operator will be denoted by $R$. These notations
    will be used throughout this paper.

       The plan of this paper goes as follows.

       \noindent
       {\bf Table of Contents}\par
\smallskip
\noindent
       {\it \quad Introduction} \par
       {\it Section 1. The existence of K\"ahler structures} \par
       {\it Section 2. A result of Hamilton and some examples} \par
       {\it Section 3. The instability of some minimal $2$-spheres \par
         \quad \quad 3.1. The twisted property of the normal bundle \par
\quad \quad 3.2. The B\"ochner-Weitzenb\"ock formula for minimal surfaces \par
References}

\head Section 1. The existence of K\"ahler structures \endhead

In this section, we shall show that, if a Riemannian $4$-manifold $(M^4, g)$ satisfies
$\frac s6 - W_+ \ge 0$, (or $\frac s6 - W_- \ge 0$) and $M^4$ is of indefinite type, then the
metric $g$ has to be K\"ahler. The main technique is to find a K\"ahler structure by using
the B\"ochner-Weitzenb\"ock formula on $2$-forms. This idea is borrowed from Bourguignon and LeBrun
(cf. [B], [L$e_1$]).

    In what follows, we always assume that $M^4$ is an oriented smooth simply-
    connected closed $4$-manifold. It is well-known that the intersection form $\omega$ on
    $H_2(M, R)$ is non-degenerate and symmetric, since $M$ is $1$-connected. By $b_2^+$ (resp. $b_2^-$)
    we mean the number of positive (resp. negative) eigenvalues of $\omega$, counting multiplicities.
    For any given metric $g$ on $M$, $s$ stands for the scalar curvature of $g$. The self-dual part (resp. anti-self-dual)
    of the Weyl tensor is denoted by $W_+$ (resp. $W_-$).

    We begin with

\proclaim{Lemma 1.1} Let $(M^4, g)$ and $b^+_2$ be as above. If $\frac s6 - W_+ \ge 0$ and
$b^+_2 \neq 0$, then $(M^4, g)$ is a K\"ahler manifold. Moreover, $\frac s6 - W_+ \ge 0$ can be expressed
as
$$
    \frac s6 - W_+ = \left[ \matrix 0 & 0 & 0 \\
0 & \frac s4 & 0 \\
0 & 0 & \frac s4 \endmatrix \right] \tag1.1
$$
relative to a suitable basis. The same conclusion holds if we replace $\{b^+_2, W_+ \}$ by
$\{ b_2^-, W_-  \}   $ above.
\endproclaim
\demo{Proof} ([L$e_1$]) Since the proof of this lemma is very short, we include it here for the sake of
completeness.

   Let $*$ be  the Hodge operator defined on $\Lambda^2(M^4)$ and let
   $$
                     H^2_+(M^4, \Bbb R) = \{ \alpha \in  \Lambda^2(M^4) \quad | \quad
                     \Delta \alpha = 0, * \alpha = \alpha \}.
   $$

   Clearly, $b^+_2 = \dim_{\Bbb R} [    H^2_+(M^4, \Bbb R)  ]$. Since $b^+_2 \neq 0$,
   we have non-zero $2$-form $\alpha \in  H^2_+(M^4, \Bbb R) $. The B\"ochner-Weitzenb\"ock
   formula tells us that $\Delta = (d + * d *)^2 = \nabla^* \nabla - 2 W + \frac s3 $,
   where $\nabla$ is the covariant derivative of $g$ (cf. [FU] p213). Hence, one sees that
   $$
        0 = \int_M \| \nabla\alpha \|^2 + 2 \int_M \langle \alpha, (\frac s6 - W_+) \alpha
        \rangle.
   $$

Since $\frac s6 - W_+ \ge 0$, both terms are non-negative and so vanish.
It follows from $\nabla\alpha = 0$ that $\alpha$ is a parallel $2$-from. Hence,
$\| \alpha \| = const$ and we may assume that $\| \alpha \|^2 = 2$.
Since $* \alpha = \alpha $, for any given point $p \in M$, there is a dual
orthonormal basis $\{ \theta_1, \theta_2, \theta_3, \theta_4 \}$ of $\Lambda^1(M^4)$ so that
$$
\alpha = \theta_1 \wedge \theta_2 + \theta_3 \wedge \theta_4
$$
where $\{ \theta_1, \theta_2, \theta_3, \theta_4 \}$ is the dual of
$\{ e_1, e_2, e_3, e_4 \}$. Using $\alpha$, one can introduce an almost complex structure $J$
by letting $Je_1 = e_2$,  $Je_3 = e_4$, $Je_2 = - e_1$ and $Je_4 = - e_3$. The condition
$\nabla\alpha = 0$ implies that $\nabla J = 0$. Hence, the Nijenhuis tensor of
$J$ vanishes. It follows from the Newlander-Nirenberg theorem that $J$ is integrable
and $(M, g)$ is a K\"ahler manifold (cf. [KN] p123-124). In this case, the expression
$W_+$ will be derived in Corollary 2.11.
\qed\enddemo

Lemma 1.1 has some interesting applications.

\proclaim{Corollary  1.1} Let $M^4$ be a smooth $4$-manifold of indefinite type
and carry a smooth metric $g$. Then $g$ satisfies $\frac s6 - W_+ \ge 0$ if and only if
$(M, g)$ is K\"ahler and $s \ge 0$.
\endproclaim

\proclaim{Theorem A} Let $M^4$ be a simply-connected closed $4$-manifold whose intersection
form is indefinite. Suppose that $M^4$ carries a smooth metric $g$
of non-negative sectional curvature and
$$
\frac s6 - W_+ \ge 0, \quad \text{ or } \quad \frac s6 - W_- \ge 0. \tag*
$$
Then $M^4$ is diffeomorphic to $S^2 \times S^2$ and $g$ is a product of metric
of $K \ge 0$.
\endproclaim
\demo{Proof} K\"ahler surfaces with non-negative holomorphic bi-sectional curvature
has been classified by Howard and Smyth (cf. [HS]). The holomorphic bi-sectional
curvature can be defined as following
$$
 K^h (X, Y) = R(X, JX, Y, JY)
$$
where $J$ is the complex structure of $(M^4, g)$ and
$$
\cases R(X, Y)Z = \nabla_X \nabla_YZ - \nabla_Y\nabla_X Z - \nabla_{[X, Y]}Z \\
       R(X, Y, Z, V) = \langle R(X, Y)V, Z \rangle. \endcases
$$
Since $(M, g)$ is K\"ahler, we know that $\nabla J = 0$ and
$$
      R(X, Y, JZ, JU) = R(X, Y, Z, U) = R(JX, JY, Z, U).
$$
The Bianchi identity and $J^2 = - Id$ implies that
$$
\align
0 & =  R(X, Y, JZ, JU) + R(X, Y, JZ, JU) + R(X, Y, JZ, JU)\\
  & =  R(X, Y, JZ, JU) + R(X, Y, JZ, JU) + R(X, Y, JZ, JU)
\endalign
$$
Now, if sectional curvature are non-negative, then
$$
\align
K^h(X, Y) & =  R(X, JX, Y, JY) \\
  & =  R(X, JY, X, JY) + R(X, Y, X, Y)  \ge 0.
\endalign
$$
The rest of proof immediately follows from [HS].
\qed\enddemo

    The Howard-Smyth's result is a special case of the generalized
    Frankel's conjecture which has been  solved by N. Mok (cf. [Mok]).

    Theorem A implies that there is no metric $g$ on $S^2 \times S^2$ which satisfies $K > 0$ and
$\frac s6 - W_- \ge 0$. The Hopf problem of the non-existence of metrics on $S^2 \times S^2$ with
positive sectional curvature  $K > 0$ remains unsolved. In the next section, we will construct
$1$-parameter family of {\it twisted metrics} on $S^2 \times S^2$ with $\frac s6 - W_+ \ge 0$.

   Lemma 1.1 and the proof of Theorem A also imply

\proclaim{Theorem B  } Let $M$ be a smooth simply-connected closed $4$-manifold which carries a
metric $g$ of non-negative curvature $K \ge 0$ and $\frac s6 - W_+ \ge 0$. Then one of following holds:
\smallskip
\noindent
(1) $M^4$ is homeomorphic to $S^4$;
\smallskip
\noindent
(2) $(M^4, g)$ is complex, K\"ahler and bi-holomorphic to $\Bbb CP^2$;
\smallskip
\noindent
(3) $M^4$ is diffeomorphic to $S^2 \times S^2$, $g$ is a product of metrics of non-negative
curvature on each factor.
\endproclaim
\demo{Proof} Since $\pi_1(M^4) = 0$, $H_2(M^4) $ has no torsion. We divide the
proof into two cases.

\medskip
\noindent
{\it Case 1.} $H_2(M^4, \Bbb Z) = 0$.

    It follows from Freedman's theorem that $M^4$ has to be homeomorphic to
    $S^4$ (cf. [L] [FU]).

\medskip
\noindent
{\it Case 2.} $H_2(M^4, \Bbb Z) \neq 0$.

    Using Lemma 1.1, one sees that $(M^4, g)$ is K\"ahler. One can now apply Howard-Smyth's theorem
    as in the proof of Theorem A.
\qed\enddemo

    The curvature condition $\frac s6 - W_+ \ge 0$ (or $\frac s6 - W_- \ge 0$) will be discussed
    in the next section.

    \head Section 2. A result of Hamilton and Some Examples \endhead

    In this section, we will give an alternative proof of Hamilton's result on the classification of
    $4$-manifolds with non-negative curvature operators. We also want to construct an $1$-family of
    {\it twisted} metrics on $S^2 \times S^2$ with $\frac s6 - W_+ \ge 0$.

       First, we begin with the Singer-Thorpe description of curvature tensor (cf. [Si-Th]). Let
       $M^n$ be $n$-dimensional Riemannian manifold with tangent space $T_p(M)$ at the point $p \in
       M$. The curvature operator $R$ at $p$ is the self-adjoint linear endomorphism:
       $$
          R: \Lambda^2(T_pM) \to \Lambda^2(T_pM)
       $$
       on the second exterior power $\Lambda^2(T_pM)$ of the tangent space defined by formula:
       $$
          \langle R(X \wedge Y), U \wedge V \rangle = \langle R(X,  Y)V,  U \rangle \tag2.1
       $$
for $X, Y, U, V \in T_pM$, where the inner product $\langle , \rangle $ are defined by the Riemannian
metric $g$ and $R$ is the usual Riemann-Christoffel curvature tensor defined by
(1.2).

     The curvature tensor $R$ of $g$ has the following decomposition:
$$
     R = \frac{s}{n(n-1)} g \oslash g + \frac{2}{n-2} \hat{R}ic \oslash g + W
$$
where $\hat{R}ic (X, Y) = Ric(X, Y) - \frac sn g$, $Ric (X, Y) =
- trace \{ Z \to R(X, Z) Y \}$ and for any bilinear symmetric form $B$, we let
$$
  (B \oslash g) (X, Y; V, W)  = \frac 12
  \{ \det \left[ \matrix  g(X, V), &  g(X, W) \\
                                                             B(Y, V), & B(Y, W) \endmatrix \right]
                                                             +
                                                             \det \left[ \matrix  B(X, V), &  B(X, W) \\
                                                             g(Y, V), & g(Y, W) \endmatrix \right]\}.
$$

   Now suppose $n = 4$ and let $M$ be given an orientation. Then the Hodge star operator $*:
   \Lambda^2(T_pM) \to \Lambda^2(T_pM)$ has two eigenvalues $\pm 1$. It is clear that $ *^2 = Id$.
   For any $\alpha \in \Lambda^2(T_pM)$, one has $ \alpha = \frac 12 (\alpha + * \alpha)
   + \frac 12 (\alpha - * \alpha)$. Hence, $\Lambda^2(T_pM) = \Lambda^2_+(T_pM) \oplus \Lambda^2_-(T_pM)$,
   where $\Lambda^2_{\pm}(T_pM) = \{ \beta \in \Lambda^2(T_pM) \quad | \quad
   * \beta = \pm \beta \}$. Singer and Thorpe show that $ * W * = W$. Therefore, the Weyl tensor decomposes
   into two parts $W_\pm$.
$$
  W \left[\matrix \Lambda^2_+(T_pM) \\ \Lambda^2_-(T_pM) \endmatrix \right]
= \left[\matrix W_+, & 0  \\ 0, & W_- \endmatrix \right]
\left[\matrix \Lambda^2_+(T_pM) \\ \Lambda^2_-(T_pM) \endmatrix \right].
$$
Using (2.2) we get
$$
  R \left[\matrix \Lambda^2_+(T_pM) \\ \Lambda^2_-(T_pM) \endmatrix \right]
= \left[\matrix  \frac{s}{12} + W_+, & \hat{R}ic  \\ [\hat{R}ic ]^t , &  \frac{s}{12} +  W_- \endmatrix \right]
\left[\matrix \Lambda^2_+(T_pM) \\ \Lambda^2_-(T_pM) \endmatrix \right], \tag2.3
$$
where  $B^t$ is the transpose of $B$.

The following observation was made in [MW] and [MM].

\proclaim{Lemma 2.1  } ([MW]) (i) If $\frac{s}{12} +  W_+ \ge 0$, then $\frac{s}{6} -  W_+ \ge 0$.
(ii) If $\frac{s}{12} +  W_- \ge 0$, then $\frac{s}{6} - W_- \ge 0$.
\endproclaim
\demo{Proof} The proof can be found in [MM] and [MW]. Since we want to use this proof later on, we
reproduce it here.

    Let $\{e_1, e_2, e_3, e_4 \}$ be a locally defined and positively oriented orthonormal frame of
    $(M^4, g)$ and set
    $$
    \aligned
& \eta_1 = e_1 \wedge e_2 + e_3 \wedge e_4  \\
& \eta_2 = e_1 \wedge e_3 - e_2 \wedge e_4  \\
& \eta_3 = e_1 \wedge e_4 + e_2 \wedge e_3
\endaligned
    $$
    Clearly, $\{ \frac{\eta_1}{2}, \frac{\eta_2}{2}, \frac{\eta_3}{2} \}$ form an orthonormal
    basis for $\Lambda^2_+(T_pM)$.

    Since $W_+$ is symmetric operator defined on $\Lambda^2_+(T_pM)$, for a given point $p$, we
    may assume
    $$
 W_+  \left[\matrix  \eta_1 \\ \eta_2 \\ \eta_3 \endmatrix \right]
= \left[\matrix  \lambda_1, & 0,  &  0 \\  0, & \lambda_2  &  0
\\  0, & 0, &  \lambda_3 \endmatrix \right]
 \left[\matrix  \eta_1, \\ \eta_2, \\ \eta_3 \endmatrix \right].
 $$
It is known that
$$
0 = tr (W)  = tr(W_+) + tr(W_-) \tag2.4
$$
and $W$ also satisfies the Bianchi identity, which implies
$$
\aligned
0 = tr(W*)  & = tr(W_+*) + tr(W_-*) \\
& = tr(W_+) - tr(W_-).
\endaligned \tag2.5
$$
Combining (2.4) and (2.5), one gets
$$
tr(W_+) = tr(W_-) = 0. \tag2.6
$$
Hence, $\lambda_1 + \lambda_2 + \lambda_3 = 0$.

(i) If $\frac{s}{12} + W_+ \ge 0$, one sees that
$$
   0 \le ( \lambda_1 + \frac{s}{12} ) + (\lambda_2 + \frac{s}{12} ) = \frac s6 - \lambda_1.
$$

    For the same reason, one gets $\frac s6 - \lambda_1 \ge 0 $, $ i = 1, 2, 3$. Thus, we proved that
    $
    \frac s6 - W_+ \ge 0.
    $

    (ii) We leave the proof to the reader. \qed \enddemo

    Using Lemma 2.1 and Theorem B, we have
\proclaim{Corollary 2.2 } (Hamilton [H]) Let $M^4$ be a smooth closed $1$-connected $4$-manifold. If $M^4$ carries a metric
of non-negative curvature operator. Then one of following holds:
\smallskip
\noindent
(1)  $M$ is homeomorphic to $S^4$;
\smallskip
\noindent
(2) $(M,g)$ is complex, K\"ahler and bi-holomorphic to $\Bbb CP^2$;
\smallskip
\noindent
(3) $M^4$ is diffeomorphic to $S^2 \times S^2$, $g$ is a product of metrics of non-negative curvature on each factor.
\endproclaim
\demo{Proof} The condition of non-negative curvature operator $R \ge 0$ implies that $\frac{s}{12} + W \ge 0$. Hence the corollary
follows from Lemma 2.1 and Theorem B immediately.
\qed\enddemo

   In case (i), R. Hamilton indeed showed that $M^4$ is diffeomorphic to $S^4$. Our proof does not yield this diffeomorphism result, but only for
   the homeomorphism part.
   Using the classification of holonomy group, Gallot and Meyer derived some relevant results (cf. [GM]).

   As we mentioned in the introduction, any closed $1$-connected smooth $4$-manifold of definite type is homotopy equivalent
   to a Riemannian manifold $(M^4, g)$ of $\frac{s}{6} - W_+ \ge 0$. In what follows, we consider $4$-manifold of indefinite type,
   for example, $S^2 \times S^2$. We will construct an $1$-parameter family of twisted metrics $g_t$ on $S^2 \times S^2$ whose volumes are equal
   to $16 \pi^2$ and $g_t$ satisfies $\frac{s}{6} - W_+ \ge 0$ for every $t$.

      It will take several steps to express these metrics explicitly. We begin with

\proclaim{Lemma 2.8} ([MW]) If $(M^4, g)$ is a K\"ahler surface, then $W_+$ can be expressed as
$$
W_+ = \left[\matrix  \frac s6, & 0,  &  0 \\  0, & \frac{-s}{12},  &  0
\\  0, & 0, &  \frac{-s}{12} \endmatrix \right] \tag2.9
$$
relative to a suitable basis. Furthermore, if $\eta_1$ is the dual of K\"ahler class, then $W_+(\eta_1) = \frac s6 \eta_1$.
\endproclaim
\demo{Proof} For any given point $p \in M^4$, we can choose a local holomorphic coordinate system $\{ z_1, z_2 \}$ such that
$$
\aligned
& Z_1 = \frac{\partial}{\partial z_1} |_p = \frac{e_1 + i e_2}{\sqrt 2} \\
& Z_1 = \frac{\partial}{\partial z_2} |_p = \frac{e_3 + i e_4}{\sqrt 2}
\endaligned
$$
where $\{e_1, e_2, e_3, e_4\}$  is an oriented orthonormal basis of $T_pM$.

   It is easy to check that
   $$
\aligned
& Z_1 \wedge \bar{Z}_1 + Z_2 \wedge \bar{Z}_2 = i (e_1 \wedge e_2 + e_3 \wedge e_4) \\
& Z_1 \wedge {Z}_2  = \frac 12 (e_1 \wedge e_3 -  e_2 \wedge e_4)+  \frac i2 (e_2 \wedge e_3 + e_1 \wedge e_4).
\endaligned
   $$
   Hence we get
   $$
\aligned
& \eta_1 = e_1 \wedge e_2 + e_3 \wedge e_4 = i (   Z_1 \wedge \bar{Z}_1 + Z_2 \wedge \bar{Z}_2        ) \\
& \eta_2 = e_1 \wedge e_3 - e_2 \wedge e_4 =    Z_1 \wedge {Z}_2 + \bar{Z}_1 \wedge \bar{Z}_2         \\
& \eta_3 = e_1 \wedge e_4 + e_2 \wedge e_3 = i (   Z_1 \wedge {Z}_2 - \bar{Z}_1\wedge \bar{Z}_2        ).
\endaligned
   $$
It is also known that $R(Z_1 \wedge {Z}_2) = R(\bar{Z}_1\wedge \bar{Z}_2) = 0$, since $(M^4, g)$ is K\"ahler (cf. [KN] p155-159). Thus,
$$
R(\eta_2) = R(\eta_3) = 0. \tag2.10
$$
Using (2.10) and the fact that $tr(W_+) = 0$, we conclude
$$
\aligned
  \langle R(\eta_1), \eta_1 \rangle & = 2 tr(R|_{\Lambda^2_+(T_pM)}) - \langle R(\eta_2), \eta_2 \rangle -
   \langle R(\eta_2), \eta_2 \rangle \\
   & = 2 tr( \frac{s}{12} + W_+) = \frac s2,
   \endaligned \tag 2.11
$$
$$
R|_{\Lambda^2_+(T_pM)} =   \left[\matrix  \frac s4 &   &  \\   & 0  &
\\   &  &  0 \endmatrix \right]
$$
and
$$
W_+ = R - \frac{s}{12} =   \left[\matrix  \frac s6 &   &  \\   & \frac{-s}{12}  &
\\   &  &  \frac{-s}{12} \endmatrix \right].
$$
This completes the proof of Lemma 2.8. \qed\enddemo

\proclaim{Corollary 2.10} ([L$e_1$], [Der]) A K\"ahler surface is anti-self-dual (i.e. $W_+ = 0$) if and only if $s= 0$.
\endproclaim

\proclaim{Corollary 2.11}  If  $(M^4, g)$  is a K\"ahler surface, then
$$
 \frac{s}{6}  - W_+ =   \left[\matrix  0 &   &  \\   & \frac{s}{4}  &
\\   &  &  \frac{s}{4} \endmatrix \right].
$$
\endproclaim
\demo{Proof} It follows from (2.3) and (2.9) that
$$
\aligned
 \frac{s}{6}  - W_+   & = \frac s6 - (R - \frac{s}{12}) = \frac s4 - R|_{\Lambda^2_+(T_pM)}  \\
   & =  \left[\matrix  0 &   &  \\   & \frac s4  &
\\   &  &  \frac s4 \endmatrix \right].
   \endaligned
$$
This completes the proof. \qed\enddemo

\medskip
\medskip
\noindent
{\bf Example 2.12.} First, we consider $S^2 \times S^2$. Let $h_0$ be the product of standard metric on each factor,
say $h_0 = \rho_1 \oplus \rho_2$. Suppose that $\varphi$ is a smooth function defined on $S^2 \times S^2$ whose mixed derivatives are
not zero. This means that if $(z_1, z_2)$ is a local coordinate system of $\Bbb CP^1 \times  \Bbb CP^1 = S^2 \times S^2$, then
$ \frac{\partial^2 \varphi }{\partial z_1 \partial \bar{z}_1} \neq 0 $ at some points of $S^2 \times S^2$.

We define $h_t = \{ (1 - \frac{t^2}{4}) \rho_1 \oplus (1 - \frac{t^2}{4})^{-1} \rho_2  \}$ and
$$
g_t = h_t + i \varepsilon \partial \bar{\partial} \varphi,
$$
where $\partial \bar{\partial} \varphi = \Sigma \frac{\partial^2 \varphi }{\partial z_\alpha \partial \bar{z}_\beta}
d z_\alpha \wedge d \bar{z}_\beta$.

It can be shown that $g_t$ is {\it not} a product metric as long as $\varepsilon \neq 0$. Since $g_t = $ is
cohomologuous to $h_t$, the K\"ahler class of $g_t$ varies as $t$ changes. However, it follows from the Stokes' theorem
that $vol(M, g_t) = vol(M, h_t) = 16 \pi^2$. Furthermore, for a fixed $\varphi $, the scalar curvature of $g_t$ is strictly
positive when $\varepsilon \to 0$. Hence, $\frac s6 - W_+ \ge 0$ holds for any $t \in [0, 1]$ and sufficiently small $\varepsilon$.

   On $\Bbb CP^2 \# k \overline{\Bbb CP}^2$, Tian and Yau constructed positive K\"ahler Einstein metrics when $3 \le k \le 8$. Using this fact
   and the same method as above, one can get an $1$-family of metrics $g_t$ of $\frac s6 - W_+ \ge 0$, which are not
   Einstein metrics.

   \head Section 3. The instability of some minimal $2$-spheres \endhead

   In this section, we shall study the instability of a minimal $2$-sphere in a $4$-manifold with positive sectional curvature.

      We are interested in this problem since there always exist some minimal $2$-spheres in $(M^n, g)$ which generate
      $\pi_2(M^4)$. This result was due to Sacks-Uhlenbeck  and improved by Meeks-Yau and Siu-Yau. Furthermore,
      using the minimal $2$-sphere theory, Siu and Yau solved the Frankel's conjecture. Micallef and Moore
      obtained a very interesting result for real manifolds (cf. [SiY], [MM]). For a given Riemannian
      $4$-manifold $M^4$, a result of Futaki indicates that the existence of the area of the area
      minimizing $2$-sphere within a given homotopy class is related to the intersection form of
      $M^4$ (cf. [Fut]). Our goal in the long run is to use the minimal $2$-sphere theory to
      study intersection forms of $4$-manifolds.

         First, we want to begin with elementary examples to illustrate the
         main idea of this section. Suppose that $g_0$ is the standard
         metric on $\Bbb CP^2$ and $S = \Bbb CP^1 \subset \Bbb CP^2$ be the
         usual totally geodesic $2$-sphere. Applying the Stokes theorem to the
         K\"ahler form $\omega$, one can show that $\Bbb CP^1$ is area-minimizing
         among its homology class. We remark that the normal bundle of $\Bbb CP^1$ in
         $\Bbb CP^2$ is non-trivial. Let us now consider another example. If the
         unit $2$-sphere $S^2$ is canonically embedded in the unit $4$-sphere $S^4$,
         then $S^2$ is totally geodesic and unstable, since the normal bundle
         of $S^2$ is parallel.

            Two examples indicate that, for minimal $2$-sphere $S$ in a Riemannian
            $4$-manifold $(M^4, g)$, the stability of $S$ is related to the twisted
            property of its normal  bundle.  The topology of the normal bundle of
            $S$ is clearly related to the intersection form of $M^4$ evaluated
            at the homology class of $S$. In fact, Kawai, Micallef and Wolfson have derived some
            interesting results in this direction. In order to prove Theorem C, we shall
            recall some useful facts in [Ka] and [MW], and we will also make
            some additional observation.

            \subhead 3.1. The twisted property of the normal bundle \endsubhead

            Our goal is to use the intrinsic curvature of the normal bundle to study
            the instability of minimal $2$-spheres. The intrinsic curvature of
            a normal bundle determines its own twisted properties.

            Let $S$ be an oriented surface minimal immersed in an oriented Riemannian
            $4$-manifold $(M^4, g)$ and $N(S)$ be the oriented normal bundle of $S$ by orientations
            of $S$ and $M^4$. We simply denote $N(S)$ by $\nu$. Let
$$
  R^\bot(X, Y) Z = \nabla^\bot_X \nabla^\bot_Y Z -  \nabla^\bot_Y \nabla^\bot_X Z -
\nabla^\bot_{[X, Y]} Z \tag3.1
$$
be the curvature tensor of normal bundle, where  $\nabla^\bot $ is the induced covariant
derivative on $N(S) \subset T(M)|_S$ associated with the induced metric:
$$
\nabla^\bot_X \sigma = (\nabla_X \sigma)^N \tag3.2
$$
and $( .)^N$ means the projection $( .)^N: T_p(M^4) \to N_p(S)$.

The  curvature $K^\bot $ of the normal
            bundle can be defined as
            $$
K^\bot = \langle R^\bot(e_1, e_2)e_4, e_3 \rangle
            $$
where $ \{ e_1, e_2 \}$ is an {\it orthonormal} basis of $T_p(S)$, $\{ e_3, e_4 \}$ is an {\it orthonormal} basis of $N_p(S)$ and
$\{ e_1, e_2, e_3, e_4 \}$ gives the positive orientation of $M^4$.

Our program is motivated by following observations.

\proclaim{Lemma 3.3} Let $S$ be an immersed minimal $2$-sphere in $(M^4, g)$. Then $K^\bot = 0$ if and
only if there is a non-zero a non-zero parallel normal cross section $E$ of $N(S)$, i.e.,
$E \neq 0$ and $\nabla^\bot E = 0$.
\endproclaim

\proclaim{Corollary 3.4} Let $S$  and $(M^4, g)$ be as in Lemma 3.3. If sectional
curvature of $(M^4, g)$ are positive and $K^\bot = 0$, then $S$ is unstable.
\endproclaim
\demo{Proof} The proof of Lemma 3.3 is straightforward (cf. Lemma 3.15 as well). Assume that Lemma 3.3
is true, we can find a global parallel normal cross section $E_3$ of $N(S)$ and
$\|E_3 \| = 1$.

It is well-known that the second variation of areas along $S$ in the direction $\delta^2 (E_3)$ can be expressed
as
$$
\aligned
\delta^2 (E_3) & = \int_S \{ \|\nabla^\bot E_3 \|^2 - K( E_3, e_1) - K(E_3, e_2) -
\| A^{E_3} \|^2   \} dA_S\\
& = - \int_S \{ K( E_3, e_1) + K(E_3, e_2) +
\| A^{E_3} \|^2   \}dA_S < 0
\endaligned
$$
where $K(X, Y)$ is the sectional curvature of the $2$-plane spanned by $X$ and $Y$
in $(M^4, g)$ and $ \| A^{\sigma} \| $ is the norm of the second fundamental form in the direction
$\sigma$ (cf. (1.18) in [SL]). \qed\enddemo

  When $K^\bot \neq 0$ for $N(S)$, Lemma 3.3 indicates that there is no global non-zero parallel
  normal cross section in $\nu = N(S)$. Therefore, we are going to pick up
  some normal cross sections, which are {\it ``half-parallel"} or {\it ``pseudo-holomorphic"}
  (cf. [Ka]).

     In order to carry out this idea, we need to define a complex structure of $N(S)$. Let $\{e_3,
     e_4 \}$ be any positively oriented orthonormal basis of $N_p(S)$, we let $Je_3 = e_4$ and
     $Je_4 = - e_3$. It is easy to check (cf. [Ka])
     $$
            \nabla^\bot_X ( J\sigma) = J( \nabla^\bot_X \sigma), \tag3.6
     $$
for all $X \in T_p(S)$ and any normal cross-section $\sigma$ in $\nu$. There is also
a natural complex structure $I$ defined on $S$ by the orientation and the induced metric of
$T(S)$.

  \definition{Definition 3.7} (1)  A local normal cross section $\sigma$ is called
{\it ``half-parallel"} or {\it ``pseudo-holomorphic"} at a point $p$ if
$$
\frac 12 (    \nabla^\bot_X  \sigma   +  \nabla^\bot_{IX }( J\sigma)          ) = 0, \quad \text{ for any  } \quad X \in T_pS. \tag3.8
$$

(2) For any real vector $X \in T_pS$,  by $\nabla^\bot_{iX} \sigma $ we mean $J( \nabla^\bot_X \sigma)$.
Further, we let
$$
\| \bar{\partial}^\bot \sigma \|^2 =
\frac 12 \|\nabla^\bot_e  \sigma   +  \nabla^\bot_{Ie }( J\sigma)  \|^2 \tag3.9
$$
where $e$ is any real unit vector of $T_pS$. Thus, $\sigma$ is holomorphic at $p$
if and only if $\bar{\partial}^\bot \sigma = 0 $.

(3)   Similarly, we let $\| \nabla^\bot  \sigma \|^2 = \| \nabla^\bot_e  \sigma \|^2
+ \| \nabla^\bot_{Ie}  \sigma \|^2$.
  \enddefinition

By the definition, we can show

\proclaim{Lemma 3.10} If $\sigma$ is a smooth normal cross section, then
$$
\int_S \| \nabla^\bot  \sigma \|^2 = \int_S \{ 2 \| \bar{\partial}^\bot \sigma \|^2
+ K^\bot \| \sigma \|^2  \}
$$
\endproclaim
\demo{Proof } Let $\theta$ be an $1$-form on $S$ as follows:
$$
\theta (X) = \langle \nabla^\bot_X(J  \sigma), \sigma \rangle.
$$
We claim that
$$
d\theta = \{K^\bot \| \sigma \|^2 + 2\| \bar{\partial}^\bot \sigma  \|^2 -
  \| \nabla^\bot  \sigma \|^2                \} d A_S
$$
where $dA_S$ is the area form of $S$. Since both sides are well-defined and independent
of the choices of local coordinates. For any given point $p \in S$, we can choose
a local orthonormal frame $\{ E_1, E_2 \}$ such that
$$
   E_1(p) = e_1, E_2(p) = e_2 \quad \text{ and } (\nabla_{E_1}E_2 )^\top(p) = 0.
$$
Hence we have $[E_1, E_2]|_p = 0$. By the definition and (3.6), we get
$$
\aligned
  d\theta(e_1, e_2)
  & = E_1 \theta(E_2) - E_2 \theta(E_1) \\
  & = E_1 \langle \nabla^\bot_{E_2}(J  \sigma), \sigma \rangle
      - E_2 \langle \nabla^\bot_{E_1}(J  \sigma), \sigma \rangle \\
      & = \langle \nabla^\bot_{E_1} \nabla^\bot_{E_2}(J  \sigma), \sigma \rangle
           - \langle \nabla^\bot_{E_2} \nabla^\bot_{E_1}(J  \sigma), \sigma \rangle \\
           & \quad + \langle  \nabla^\bot_{E_2}(J  \sigma), \nabla^\bot_{E_1} \sigma \rangle
             -  \langle  \nabla^\bot_{E_1}(J  \sigma), \nabla^\bot_{E_2} \sigma \rangle \\
             &  = K^\bot(p) \| \sigma \|^2 +
             2 \langle  \nabla^\bot_{E_2}(J  \sigma), \nabla^\bot_{E_1} \sigma \rangle \\
             & = K^\bot(p) \| \sigma \|^2 + 2\| \bar{\partial}^\bot \sigma  \|^2 -
  \| \nabla^\bot  \sigma \|^2.
  \endaligned
$$
This completes the proof of (3.12) and Lemma 3.10.
\qed\enddemo

We can now rewrite the second variational formula:

\proclaim{Theorem 3.13} Let $S$ be an immersed minimal surface in $(M^4, g)$, and let
$K^\bot$ be the curvature of the normal bundle $\nu$ of $S$. Then the second variation of
$S$ in the direction $\sigma$ is
$$
\delta^2 (\sigma) = \int_S \{2 \|\bar{\partial}^\bot \sigma \|^2+ [ K^\bot - K( \sigma , e) - K(\sigma, Ie)]\|\sigma \|^2
 -
\| A^{\sigma} \|^2   \} dA_S \tag3.14
$$
where $e$ is any unit tangent vector in $T_pS$.
\endproclaim
\demo{Proof} The original second variational formula is given in [SL]:
$$
\delta^2 (\sigma) = \int_S \{ \|\nabla^\bot_e \sigma \|^2 + \|\nabla^\bot_{Ie} \sigma \|^2- K( \sigma , e) - K(\sigma, Ie) -
\| A^{\sigma} \|^2   \} dA_S.
$$
It follows form Definition 3.7 (3) and (3.11) that
$$
\int_S \{ \|\nabla^\bot_e \sigma \|^2 + \|\nabla^\bot_{Ie} \sigma \|^2\} dA_S = \int_S \{2 \|\bar{\partial}^\bot \sigma \|^2
+ K^\bot \|\sigma \|^2\}dA_S.
$$
This completes the proof.
\qed\enddemo

\bigskip
\noindent
{\bf Remark:} It is well-known that, if we let
$$
  c_1 (\nu) = \frac{1}{2\pi}\int_S   K^\bot dA_S
$$
then $c_1(\nu)$ is equal to the Euler number of $\nu$. We may always assume that $c_1 (\nu) \ge 0$ appropriately
(cf. [Ka]). When $S$ is an immersed minimal sphere, the existence of a holomorphic normal cross section is given
by the Riemann-Roch Theorem. Furthermore, if $\nu$ is topologically trivial, then
$$
\int_S   K^\bot dA_S = 0.
$$
However, when $\sigma$ is holomorphic, we have
$$
 \int_S   K^\bot  \|\sigma \|^2  dA_S = \int_S   \|\nabla^\bot \sigma \|^2 dA_S \ge 0 \quad \text{ (by Lemma 3.10). }
$$
The equality holds if and only if $\sigma$ is parallel. In fact, $\|\sigma \|$ is closely related to
$K^\bot$ as follows:

\proclaim{Lemma 3.15} Let $\sigma$ be a local non-vanishing holomorphic cross section. Then
$$
  K^\bot = - \Delta_S [\log ( \|\sigma \|^2  ) ].
$$
Furthermore, when $S$ is a $2$-sphere, $K^\bot(p) = 0$ for all $p \in S$ if and only if
there is a global non-zero  holomorphic cross section $\sigma$ with $\|\sigma \| =$ constant.
\endproclaim
\demo{Proof} Since both sides of (3.16) are independent of the choice of $S$. For any given
$p \in S$, we take geodesic coordinates $\{x, y\}$ of $(S, g|_S)$ at $p$.
Let $e_1 = \frac{ \partial }{\partial x}|_p$  and  $e_2 = \frac{ \partial }{\partial y}|_p$.

Since $\sigma$ is holomorphic, one sees that
$$
\aligned
K^\bot \|\sigma \|^2(p) & =  \langle R^\bot(\frac{ \partial }{\partial x}, \frac{ \partial }{\partial y}) J\sigma,
\sigma \rangle \\
& = \langle \nabla^\bot_{e_1} \nabla^\bot_{e_2}(J  \sigma) -
\nabla^\bot_{e_2} \nabla^\bot_{e_1}(J  \sigma), \sigma \rangle \\
& = \langle - \nabla^\bot_{e_1} \nabla^\bot_{e_1}  \sigma -
\nabla^\bot_{e_2} \nabla^\bot_{e_2} \sigma, \sigma \rangle \\
& = - \Delta_S [\|\sigma \|^2] +  \|\nabla \sigma \|^2
\endaligned
$$
Now, if $\|\sigma \| =$ constant, then $K^\bot = 0$. Conversely, when $K^\bot (p) = 0$ for all $p \in S$, there always
is a global non-zero holomorphic cross section $\sigma$, which is guaranteed by the Riemann-Roch Theorem. One can also show that
$\sigma$ never vanishes, since $c_1( \nu) = 0$. The condition $K^\bot = 0$ implies that
$\log \|\sigma \|$ is a harmonic function on $S$. Hence, $\|\sigma \| =$ constant.
\qed\enddemo

We have seen that it is complicated to figure out how large $K^\bot \|\sigma \|^2$ is. In the next subsection, we want to study
$K^\bot$ extrinsically.

\subhead 3.2. The Weitzenb\"ock formula for minimal surfaces \endsubhead

In what follows, we always let $ \eta = e_1 \wedge e_2 +  e_3 \wedge e_4$, and
$\{ e_1, e_2, e_3, e_4 \}$ be  an orthonormal basis of $(M^4, g)$ with
$\{e_1, e_2 \} \subset TS$, $\{e_3, e_4 \} \subset N(S) = \nu$. We also define
$$
\| A \wedge A \|^2 = \|A^{e_3}(e_1) + A^{Je_3}(e_2)\|^2 +  \|A^{e_3}(e_2) - A^{Je_3}(e_1)\|^2. \tag3.17
$$

Now we can give  another version of the second variational formula, or the B\"ochner-Weitzenb\"ock formula
for minimal surfaces.

\proclaim{Theorem 3.18} Let $S$ be a minimal immersed minimal surface in $(M^4, g)$ with the normal bundle
$\nu$ as above. If $\sigma$ is a normal cross section, then the averaging second variations is
$$
\delta^2(\sigma) + \delta^2(J \sigma) = \int_S \{4 \|\bar{\partial}^\bot \sigma \|^2
- [\langle (\frac s6 - W_+  )\eta, \eta \rangle + \|A \wedge A \|^2 ]\| \sigma \|^2
    \} dA_S.
$$
Furthermore, $\|A \wedge A \|^2$ vanishes if and only if $S$ is totally geodesic.
\endproclaim
\demo{Proof} Using the fact that $\nabla_X \sigma = ( \nabla_X \sigma )^\bot + (\nabla_X \sigma)^\top  $,
one can verify that
$$
K^\bot = \langle R(e_1, e_2) e_4, e_3 \rangle - \langle A^{e_3}(e_1), A^{Je_3}(e_2) \rangle
+ \langle A^{Je_3}(e_1), A^{e_3}(e_2) \rangle. \tag3.20
$$
Let $Ric^\bot (\sigma) = \langle R(e_1, \sigma) \sigma, e_1 \rangle + \langle R(e_2, \sigma) \sigma, e_2 \rangle$.
It is also true that
$$
\aligned
& Ric^\bot (\sigma) + Ric^\bot (J \sigma) \\
= & \{ \langle R(e_1, e_3) e_3, e_1 \rangle +  \langle R(e_1, e_4) e_4, e_1 \rangle \\
  & + \langle R(e_2, e_3) e_3, e_2 \rangle +  \langle R(e_2, e_4) e_4, e_2 \rangle\} \|\sigma \|^2  \\
  = & \{ Ric^\bot (e_3) + Ric^\bot (e_4) \}\|\sigma \|^2.
\endaligned
\tag3.21
$$
Using the Bianchi identity and (2.3), we have
$$
\aligned
 & -2 \langle R(e_1, e_2) e_4, e_3 \rangle +   Ric^\bot (e_3) + Ric^\bot (e_4)   \\
= &  \langle R(e_1 \wedge e_3  - e_2 \wedge e_4), e_1 \wedge e_3  - e_2 \wedge e_4 \rangle \\
& + \langle R(e_1 \wedge e_4  + e_2 \wedge e_3), e_1 \wedge e_4  + e_2 \wedge e_3 \rangle
 \\
  = &  \langle (\frac{s}{12} - W_+)(e_1 \wedge e_3  - e_2 \wedge e_4), e_1 \wedge e_3  - e_2 \wedge e_4 \rangle \\
  & +
\langle (\frac{s}{12} - W_+)(e_1 \wedge e_4  + e_2 \wedge e_3), e_1 \wedge e_4  + e_2 \wedge e_3 \rangle
 \\
 = & \frac{4s}{12} - 2 \{tr(W_+) - \frac 12
 \langle W_+ (e_1 \wedge e_2 + e_3 \wedge e_4), e_1 \wedge e_2 + e_3 \wedge e_4  \rangle
 \} \\
 = & \langle (\frac s6 - W_+) (e_1 \wedge e_2 + e_3 \wedge e_4), e_1 \wedge e_2 + e_3 \wedge e_4  \rangle
\endaligned
\tag3.22
$$
where we used $tr(W_+) = 0$, $\| e_1 \wedge e_3  - e_2 \wedge e_4\|^2 = 2$ and etc.

Obviously one gets
$$
\|A^{\sigma}\|^2  + \|A^{J\sigma}\|^2  = \{\|A^{e_3}\|^2  + \|A^{e_4}\|^2  \}\|\sigma \|^2
\tag3.23
$$
and
$$
\|A \wedge A \|^2 = \|A^{e_3}\|^2  + \|A^{e_4}\|^2 +
2 \langle A^{e_3}(e_1), A^{Je_3}(e_2) \rangle -2  \langle A^{Je_3}(e_1), A^{e_3}(e_2) \rangle.
\tag3.24
$$
Using (3.17) and (3.20)-(3.24), Theorem 3.13 and the fact that
$$
\cases
\| \nabla^\bot \sigma \| =  \| \nabla^\bot (J\sigma) \| \\
 \quad \|  \sigma \| =  \| J\sigma \|
 \endcases
$$
we finally get
$$
\aligned
& \delta^2(\sigma) + \delta^2(J \sigma)\\
 = & \int_S \{ 4 \|\bar{\partial}^\bot \sigma \|^2 + [
 2 K^\bot - Ric^\bot (e_3) - Ric^\bot (e_4)   ]\|\sigma \|^2 - \|A^{\sigma}\|^2  - \|A^{J\sigma}\|^2
\} dA_S \\
 = & \int_S \{4 \|\bar{\partial}^\bot \sigma \|^2
- [\langle (\frac s6 - W_+  )\eta, \eta \rangle + \|A \wedge A \|^2 ]\| \sigma \|^2
    \} dA_S.
\endaligned
$$

  Furthermore, if $ \|A \wedge A \| = 0$, then
  $$
  A^{Je_3}(e_1)   = A^{e_3}(e_2) \quad \text{ and } \quad   A^{Je_3}(e_2)   =  - A^{e_3}(e_1).  \tag3.25
  $$
  Since $S$ is minimal, we also know that $( \nabla_{e_1}e_1  + \nabla_{e_2}e_2) ^\bot = 0 $.

  After properly choosing $\{ e_1, e_2, e_3, e_4  \}$, we may assume that
  $$
(  \nabla_{e_1}e_1) ^\bot =  - (  \nabla_{e_2}e_2) ^\bot =  \lambda e_4 \text{ and  }  (  \nabla_{e_1}e_2) ^\bot = 0  \tag3.26
  $$
  where $\lambda$ is a real number. Thus, we get
  $$
0 = \langle \pm \lambda e_4, e_3 \rangle =  \langle (  \nabla_{e_j}e_k) ^\bot , e_3 \rangle =
\langle A^{e_3}  (e_k) ^\bot , e_j \rangle
  $$
  for all $j, k = 1, 2$. It follows that $A^{e_3} = 0$. Using  $A^{e_3} = 0$ and (3.25), we conclude that
$A^{Je_3} = 0$. Hence, $A = 0$ and $S$ is totally geodesic, when $\|A \wedge A\| = 0$.
\qed\enddemo

Now we are ready to show

\proclaim{Theorem C} Let $S$ be an immersed minimal $2$-sphere in $(S^2 \times S^2, g)$, which is homotopic to
the first factor of $S^2 \times S^2$. Suppose that
there is an open set $U$ containing $S$ in which the metric $g$ satisfies
$ \frac s6 - W_+ \ge 0$ and $K > 0$, then $S$ is unstable.
\endproclaim
\demo{Proof} Since $S$ is homotopic to the first factor of $S^2 \times S^2$, the normal
bundle $\nu$ of $S$ is topologically trivial. It follows from Riemann-Roch theorem that there is
always a global non-zero holomorphic cross-section $\sigma$ (cf. [Ka] p1072). Hence, it
is also true that $\bar{\partial}^\bot(J \sigma) = 0$ because of $\nabla^\bot J = 0$. Using
Theorem 3.18 and $ \frac s6 - W_+ \ge 0$, we get
$$
 \delta^2(\sigma) + \delta^2(J \sigma) \leq 0. \tag3.27
$$
We claim that
$$
 \delta^2(\sigma) + \delta^2(J \sigma) < 0 \tag3.28
$$
by using the K\"ahler property and Bianchi identity.

Suppose that (3.28) were not true, we will derive a contraction. The number of zeros of $\sigma$ is equal to $c_1(\nu) = 0$. Hence,
$\sigma$ never vanishes since $c_1(\nu) = 0$.

By Theorem 3.18 we see that if $\delta^2(\sigma) + \delta^2(J \sigma) = 0$ then $\| A \wedge A \| = 0$. It follows
that $S$ is totally geodesic. Therefore, $\nabla_X^\bot E = \nabla_X E$ for any normal cross-section. It follows
that
$$
     K^\bot = - \langle R^\bot(e, Ie) e_3, e_4 \rangle = - \langle R(e, Ie) e_3, e_4 \rangle. \tag3.28
$$

Using Bianchi identity and (3.28), we further get
$$
K^\bot = K^h(e_1, e_3) = R(e_1, e_3, e_1, e_3) + R(e_1, e_4, e_1, e_4) = K(e_1, e_3) + K(e_1, e_4) \tag3.29
$$

 Thus, (3.29) and our assumption $K > 0$ imply that
$$
  K^\bot > 0 \tag3.30
$$
and
$$
c_1 (\nu) = \frac{1}{2\pi} \int_S K^\bot dA_S > 0
$$
which is impossible, since the normal bundle $\nu$ is trivial and $c_1(\nu) = 0$.

Thus, (3.28) is verified . This completes the proof of Theorem C.
\qed\enddemo

\bigskip
\noindent
{\bf Remark:} In the proof of Theorem C, one can actually show that $S$ has index $2$ as follows:
$$
\delta(\zeta_1, \zeta_2) =
\int_S \{\langle \nabla^\bot\zeta_1,   \nabla^\bot\zeta_2 \rangle -
\sum_{k=1  }^2 [ \langle R(e_k, \zeta_1)\zeta_2, e_k  \rangle
 \langle A^{\zeta_1}(e_k), e_k  \rangle \langle A^{\zeta_1}(e_k), e_k  \rangle ]     \} dA_S.
$$

Since it has been proved above that
$$
\delta^2(\sigma) + \delta^2(J \sigma) < 0, \tag3.28
$$
we may assume that $\delta^2(\sigma) < 0$. When $\delta(\sigma, J\sigma) \le 0$,
one notices that
$$
\delta^2(\sigma + J\sigma) = \delta^2(\sigma) + \delta^2(J \sigma) + 2\delta(\sigma, J\sigma) <0.
\tag3.29
$$
In the case that $\delta(\sigma, J\sigma) > 0$, one can also show that
$\delta^2(\sigma - J\sigma) < 0 $. In either case, we have two linearly independent cross sections $\sigma$ and $\sigma \pm J \sigma$ such that $\delta^2(\sigma) < 0 $ and $\delta^2(\sigma \pm J\sigma) < 0 $.
Hence, index(S) $\ge 2$.

\bigskip
\noindent
{\bf Acknowledgement:} The author would like to express his deep gratitude to Professors E. Calabi, C. Croke and H. Gluck
for their interest and support in this work.  He is also very grateful to M. Micallef and D. Moore for helpful
conversations; their work is used in this paper extensively. The author wanted to thank C. LeBrun and J-P Sha for bring [H],
[Ka] and other references to his attention.

This paper was written in 1993. Over the years, the author received many requests for this preprint. The author is very grateful to
Professor Walter Wei and Professor En-Bing Lin for their efforts to publish this paper in ``International Journal of Geometry and Topology".
He also appreciated the referee's careful proof-reading.

\bigskip
\noindent
{\bf Added in Proof:} The current version of the paper was first written when I was at Institute of Advanced
Study, Princeton 1989-1990. I posted the paper on www.arxiv.org as math.DG/0701742 in January 2007.
I am grateful to Matt Gursky and  Johann Davidov for their interest in my paper and for bringing [Gu],
[GuL] and [AD] to my attention. There have been other significant progress made in the field after 1990.
Since the proof of theorems in this paper does not need to use the results
after 1993, I did not cite all recent results in the reference section.


\Refs
 \nofrills{References}
 \widestnumber\key{APS12}

 \vskip3mm

 \ref
 \key A
 \by Aminov, Ju. A
 \pages 359-375
 \paper On the instability of minimal surfaces in an $n$-dimension Riemannian space of positive curvature
 \jour Math. USSR Sbornik
 \vol {\bf 29} no  3
 \yr 1976
 \endref


 \ref
 \key AD
 \by Apostolov, V. and Davidov, J.
 \pages 438-451
 \paper Hermitian surfaces and isotropic curvature
 \jour Illinois J. Math.
 \vol {\bf 44} no  2
 \yr 2000
 \endref




 \ref
 \key  Be
 \by Besse, A. L
 \pages
 \book Einstein Manifolds
 \publ Ergebnisse der Mathematik und ihrer Grenzgebiete
 (3) [Results in Mathematics and Related Areas (3)], 10. Springer-Verlag, Berlin, 1987. xii+510 pp. ISBN: 3-540-15279-2
 \endref

 \ref
 \key Bo
 \by Bourguignon, J. P.
 \pages 263-286
 \paper Les varieties  de dimension $4$ a signature non nulle dont la courbure est harmonique sont d'Einstein.  (French)
\jour Invent. Math
 \vol {\bf 63}
 \yr 1981
 \endref

 \ref
 \key CG
 \by Cheeger, J. and Gromoll, D
 \pages 413-443
 \paper On the structure of complete manifolds of nonnegative curvature
 \jour Ann. of Math
 \vol  {\bf 96} no. 3
 \yr 1972
 \endref

 \ref
 \key Ch
 \by  Cheeger, J
 \pages 623-628
 \paper Some examples of manifolds of nonnegative curvature
 \jour J. Diff. Geom.
 \vol {\bf 8}
 \yr 1973
 \endref


 \ref
 \key Der
 \by  Derdzinnski, D
 \pages 405- 433
 \paper Self-dual K\"ahler manifolds and Einstein manifolds of dimension 4
 \jour Comp. Math
 \vol {\bf 49}
 \yr 1983
 \endref

 \ref
 \key DF
 \by  Donaldson, S. and Friedman, R
 \pages 197-239
 \paper Connected sum of self-dual manifolds and deformations of singular space
 \jour Nonlinearity
 \vol {\bf 2}
 \yr 1989
 \endref



 \ref
 \key FU
 \by  Freed, D. and Uhlenbeck, K
 \book Instantons and four-manifolds
 \publ Mathematical Sciences Research Institute Publications, 1. Springer-Verlag, New York, 1984. viii+232 pp. ISBN: 0-387-96036-8
 \endref


 \ref
 \key Fut
 \by  Futaki, A
 \pages 291-293
 \paper Nonexistence of minimizing harmonic maps from $2$-spheres
 \jour Proc. Japan Acad.
 \vol {\bf 56} no. 6
 \yr 1980
 \endref


\ref
 \key GM
 \by  Gallot, S. and Meyer, D
 \pages 259-284
 \paper Operateur de courbure et laplacien des formes differentielles d'une varieties riemannienne.  (French)
 \paperinfo  J. Math. Pures Appl. (9)  54  (1975), no. 3
\endref

\ref
 \key G$r_1$
\by  Gromov, M
 \pages 153-179
 \paper Curvature, diameter and Betti numbers
 \jour Math. Helv.
 \vol {\bf 56}
 \yr 1981
 \endref

\ref
 \key G$r_2$
 \by  Gromov, M
 \pages 307-347
 \paper Pseudo holomorphic curves in sympletic manifolds
 \jour Invent. Math.
 \vol {\bf 82}
 \yr 1985
\endref

\ref
 \key Gu
 \by  Gursky, M.
 \pages 417-431
 \paper Four-manifolds with
$\delta W\sp +=0$ and Einstein constants of the sphere
 \jour  Math. Ann.
 \vol {\bf 318}
 \yr 2000
\endref



\ref
 \key GuL
 \by  Gursky, M. and LeBrun, C.
 \pages 315-328
 \paper On Einstein manifolds of
positive sectional curvature
 \jour  Ann. Global Anal. Geom.
 \vol {\bf 17}
 \yr 1999
\endref


 \ref
 \key HS
 \by  Howard, A. and Smyth, B
 \pages 491-501
 \paper K\"ahler surfaces of nonnegative curvature
 \jour J. Diff. Geom.
 \yr 1971
 \vol {\bf 5}
 \endref


 \ref
 \key Ka
 \by  Kawai, S
 \pages 1067-1075
 \paper On the instability of a minimal surface in a $4$-manifold whose curvature lies in the interval $(\frac 14,\,1]$.
 \jour Publ. Res. Inst. Math. Sci of Kyoto Univ.
 \yr 1982
 \vol {\bf 18} no. 3
 \endref



 \ref
 \key KN
 \by Kobayashi, S. and Nomizu, K
 \book Foundations of Differential Geometry, vol II
 \publ Interscience Publishers, 1969
 \endref

 \ref
 \key L
 \by  Lawson, H. B
\book The theory of gauge fields in four dimensions
\publ CBMS Regional Conference Series in Mathematics, 58. Published for the Conference Board of the Mathematical Sciences,
Washington, DC; by the American Mathematical Society, Providence, RI, 1985. vii+101 pp. ISBN: 0-8218-0708-0
 \endref

\ref
 \key L$e_1$
 \by  LeBrun, C
 \pages 637-640
 \paper On the topology of self-dual $4$-manifolds
 \jour Proc. of Amer. Math. Soc
 \yr 1986
 \vol {\bf 98} no. 2
 \endref




\ref
 \key L$e_2$
 \by  LeBrun, C
 \pages 223-253
 \paper Explicit self-dual metrics on $\Bbb CP_2 \# ... \# \Bbb CP_2$
 \jour J. Differential Geom.
 \yr 1991
 \vol {\bf 34} no. 1
 \endref



 \ref
 \key MM
 \by Micallef, M. and Moore J. D
 \pages 199-227
 \paper Minimal two-spheres and the topology of manifolds with positive curvature on
 isotropic two-planes
 \jour Ann of Math
 \vol {\bf 127}
 \yr 1988
 \endref

 \ref
 \key MW
 \by Micallef, M. and  Wolfson, J.
 \pages 245-267
 \paper The second variation of area of minimal surfaces in four-manifolds
 \jour Math. Ann.
 \vol {\bf 295}
 \yr 1993
 \endref


 \ref
 \key Mok
 \by  Mok, N.
 \paper The uniformization theorem for compact K\"ahler manifolds of nonnegative holomorphic bisectional curvature
 \jour  J. Diff . Geom.
 \vol {\bf 27}, no. 2
  \pages 179-214
 \yr 1988
 \endref


\ref
 \key P
 \by  Poon, Y. S.
 \paper Compact self-dual manifolds with positive scalar curvature
 \jour  J. Diff . Geom.
 \vol {\bf 24}
  \pages 97-132
 \yr 1986
 \endref




 \ref
 \key SU
 \by Sacks, J and Uhlenback, K.
 \pages 1-24
 \paper The existence of minimal immersions of $2$-spheres
 \jour Ann. of Math
 \vol {\bf 113}
 \yr 1981
 \endref




\ref
 \key ScY
 \by Schoen, R.; Yau, Shing Tung
 \pages 127-142
 \paper Existence of incompressible minimal surfaces
 and the topology of three-dimensional manifolds with nonnegative scalar curvature
 \jour Ann. of Math.
 \vol {\bf 110}, no. 1
 \yr 1979
 \endref



\ref
 \key SY
 \by Sha, Ji-Ping and Yang, DaGang
 \pages 127-137
 \paper Positive Ricci curvature on the connected sums of $S^n \times S^m$.
 \jour J. Diff. Geom.
 \vol {\bf 33}, no. 1
 \yr 1991
 \endref





 \ref
 \key SiY
 \by Siu, Y-T and Yau, S. T.
 \paper Compact K\"ahler manifolds of positive bisectional curvature
 \yr 1980
 \jour Invent. Math.
  \vol {\bf 59},
 \pages 189-204
 \endref


\ref
 \key Si-Th
 \by Singer, I. M.; Thorpe, John A.
\book Lecture notes on elementary topology and geometry
\publ Scott, Foresman and Co., Glenview, Ill. 1967 v+214 pp
 \endref



\ref
 \key SJ
 \by Simon, J.
 \pages 62-105
 \paper Minimal varieties in Riemannian manifolds
 \jour Ann. of Math.
 \vol {\bf 88}
 \yr 1968
 \endref




 \ref
 \key SL
 \by Simon, Leon
 \pages 3-52
 \paper Survey lectures on minimal submanifolds
\paperinfo Ann. of Math. Stud., 103, Princeton Univ. Press, Princeton, NJ, 1983
\endref




\ref
 \key TY
 \by Tian, Gang and Yau, S. T.
 \pages 175-203
 \paper  K\"ahler-Einstein metrics on complex surfaces with $C_1 > 0$
 \jour Commun. Math. Phys.
 \vol {\bf 112}
 \yr 1987
 \endref


 \endRefs
 \enddocument

 \end